

\baselineskip=14pt
\parskip=10pt

\magnification=\magstephalf

\def\1{{\overline{1}}}
\def\2{{\overline{2}}}
\parindent=0pt
\overfullrule=0in

\def\frac#1#2{{#1 \over #2}}
\centerline
{\bf  
Using Symbolic Computation to Analyze some Children's Board Games
}
\bigskip
\centerline
{\it Shalosh B. EKHAD and Doron ZEILBERGER}
\bigskip

{\bf Abstract.} In a delightful article that recently appeared in Mathematics Magazine, David and Lori Mccune
analyze the  board game ``Count Your Chickens!", recommended to children three and up.
Alas, they use the advanced  theory of Markov chains, that presupposes a knowledge of linear algebra,
that few three-years-olds are likely to understand. Here we present a much simpler, more intuitive, approach,
that while unlikely to be understood by three-year-olds, will probably be understood by
a smart 14-year-old.
Moreover, our approach accomplishes much more, and is more efficient. It uses symbolic, rather than numeric computation.
The article is accompanied by a general Maple package, {\tt CountChickens.txt}, that can handle, in a few seconds,
{\it any} such game, not just this particular one. It is also accompanied by an even more general
Maple package  {\tt UmbralMarkov.txt} that handles any ``weighted" (discrete time) Markov chain with any number
of absorbing states.

{\bf The Maple packages.} This article is accompanied by two Maple packages 
{\tt CountChickens.txt}  and {\tt UmbralMarkov.txt} 
that can be obtained, 
along with numerous input and output files, from the front of this article

{\tt http://www.math.rutgers.edu/\~{}zeilberg/mamarim/mamarimhtml/board.html} \quad .

{\bf The Count Your Chickens! board game}

The board game {\it Snakes and Ladders} (that became ``Chutes and Ladders" in the USA, since snakes are too scary)
is too stressful for the gentle soul of a  typical three-year-old,
because it has a winner, and hence a loser.
Even {\it CandyLand} that involves
picking colored cards, rather than spinning a spinner, is not recommended, since it suffers from the same problem
and three-years-old (and not only) hate to lose, making them cry.
Hence game inventor Peggy Brown came up with a fun, stress-free, `cooperative' game [B] for kids, 
where there is only one team and `everyone wins together and loses together' (so it is really a solitaire game) called
``Count Your Chickens!" manufactured and marketed by the Peaceable Kingdom toy company.

In a delightful article that appeared recently in Mathematical Magazine, the mathematical couple David and Lori Mccune,
who play this game with their young children, use the sophisticated theory of Markov Chains, that entails a
knowledge of matrices - and matrix inverses - to compute the probability of winning, as well as the expected number of chicks at
the end. They got $0.6410$ for the former  and $39.22$ for the latter.
Our, simpler, faster, and more efficient approach agrees with their probability, but gave the more precise value of
$0.6410373996231\dots$, and got a slightly higher value for the expected number of chicks, namely
$39.32230439142343\dots$. [MM]'s stated value of $39.22$ rather than the correct $39.322$ is probably a misprint.

One of us (DZ) wrote a  Maple package, {\tt CountChickens.txt}, mentioned above, that 
enabled the other author (SBE) to find these quantities for
{\it any} such kind of board game, and go far beyond mere probability of winning and expected number of chicks.
It uses {\it symbol crunching} rather than {\it number crunching}, and
has many fewer `states', making the computations extremely fast.

But let us first define an `abstract' Count Your Chickens! game.

Let $N$ and $K$ be two positive integers. The game consists of

$\bullet$ A board with $N+1$ squares where the $1^{st}$ location is the starting place of Mama Chicken and
$N+1$ is the terminal square. Each square is either empty or labeled with one of  $K$ animals.

$\bullet$ a spinner with $K+1$  choices, all equally likely, labeled by the $K$  animals, plus an extra one called the Fox.

$\bullet$ a subset of $\{2, \dots , N+1 \}$ called the set of {\it blue squares}.

The rules are as follows. Mama Chicken starts out at location $1$.
At every turn, the player spins the spinner. If it is a Fox, then
you lose a chick (if you currently have no chicks, then nothing happens) and stay where you are.
Otherwise you go to the {\bf next} location labeled by the animal that you got.
The three-year-old counts the number of squares moved and collects that number of chicks.
If the new location is a blue square, then you get an extra chick.

Sooner or later, with probability $1$, you would get to the terminal square, $N+1$ that is labeled by all the $K$ animals.

You  win the game if you have at least $N$ chicks, and  otherwise you lose.

In the simplified example of [MM], N=8, K=2, the board is
$$
[START,EMPTY, SHEEP, COW, EMPTY, COW, EMPTY, SHEEP , \{COW,SHEEP\}]
$$
and the set of blue squares is $\{3,6\}$.

In the actual  game[B], $K=5$ and $N=40$. The board is as follows
$$
[0,0,S,P,T,C,D,P,C,D,S,T,0,C,P,0,0,0,T,0,T,D,S,C,D,P,T,0,S,C,0,0,T,P,S,D,0,S,C,P, 
$$
$$
\{C,D,P,S,T\}] \quad,
$$
(where C:=Cow, D:=Dog,  P:=Pig, S:=Sheep,   T:=Tractor, and $0$ indicates an empty square)
and the set of blue locations is
$$
 \left\{ 5,9,23,36,40 \right\} \quad .
$$

In [MM] the game is modeled as a Markov chain with a huge number of states, 
each of the form (location, Current Number of Chicks), essentially $O(N^2)$. For the
problem of just computing the probability of winning (for $N=40$), they manage to reduce it to
$163$ states, but for the harder problem of computing the expected number of chicks at the end, they needed $668$ states,
and the matrices were huge.

Our approach  also makes use of Markov chains, but we don't need any of the standard theory, and we never mention the word `matrix'.
Also our number of states is $O(N)$  (obviously the EMPTY squares can be ignored). We use Gian-Carlo Rota's seminal idea
of an {\it umbral operator}.

Let $f_i(t)$ be the probability generating function of landing at square $i$, where the coefficient of $t^j$ is
the probability that you currently have $j$ chicks. To indicate the fact that it is currently at location $i$ we will denote it
by $s^i f_i(t)$. If you got a Fox this becomes $s^i f_i(t)/t$ (followed by replacing $t^{-1}$ by $1$, if necessary).
Otherwise, Mama Chicken goes to a new location, let's call it $j$, and the new state becomes
$s^j f_i(t) t^{j-i}$ if $j$ is not a blue square, and  $s^j f_i(t) t^{j-i+1}$ if it is.
If we get a power of $t$ larger  than $N$, we replace it by $t^N$.

This introduces an `evolution operation' that we call the {\it pre-umbra}.

In the simplified game used in [MM], (whose board was given above), we have
$$
s^1 \rightarrow \frac{1}{3} (s^1 + s^3\,t^{3-1+1} + s^4\,t^{4-1} ) \, = \, \frac{1}{3} (s+ s^3 t^3+ s^4t^3 ) \quad ,
$$
$$
F(t) s^3 \rightarrow \frac{F(t)}{3} (\frac{s^3}{t}+ s^4\,t^{4-3} + s^6\,t^{6-3+1} )\, = \, \frac{F(t)}{3} (\frac{s^3}{t}+ s^4 t+ s^6t^4) \quad ,
$$
$$
F(t) s^4 \rightarrow \frac{F(t)}{3} (s^4/t+ s^6\,t^{6-4+1} + s^8\,t^{8-4} )\, = \, \frac{F(t)}{3} (\frac{s^4}{t}+ s^6 t^3+ s^8 t^4) \quad ,
$$
$$
F(t) s^6 \rightarrow \frac{F(t)}{3} (s^6/t+ s^8\,t^{8-6} + s^9\,t^{9-6} ) \, = \, \frac{F(t)}{3} ( \frac{s^6}{t} + s^8 t^2+ s^9 t^3) \quad ,
$$
$$
F(t) s^8 \rightarrow \frac{F(t)}{3} (s^8/t+ s^9\,t^{9-8} + s^9\,t^{9-8} ) \, = \, \frac{F(t)}{3} ( \frac{s^4}{t}+ s^9 t+ s^9 t^2) \quad ,
$$
$$
F(t) s^9 \rightarrow F(t)s^9 \quad .
$$
(since $9$ is an {\it absorbing state}).

These operations must be followed by a ``clean-up" operation. Replacing $t^{-1}$ by $1$ (you can't have a negative number of chicks),
and replacing $t^9, t^{10}, \dots$ by $t^8$.

This is the {\it pre-umbra}, defined on every {\bf monomial} $s^i$, let's call it $T$.
If we have a {\bf polynomial} in $s$ (and of course $t$), we  extend  it by {\bf linearity}.
(Recall that every {\bf poly}nomial is a linear combination of {\bf mono}mials).
We call this linear extension the {\it umbra} and also denote it by $T$.

It is readily seen that applying this operator, starting with the initial state $s^0$, describes the
`evolution' of the process. 

While, in principle, the game can last forever (if you are really unlucky, you may keep getting foxes), life is finite, so we decide that
we are playing at most $M$ rounds, and make $M$ large enough so that the probability of lasting longer than $M$ rounds is negligible.

The probability generating function after $1$ round is $T(s^0)=\frac{1}{3} (1+ s^2 t^3+ s^3t^3 )$.
After two rounds is $T^2(s^0)$, etc.. Sooner or later we will encounter $s^9$ (in general $s^{N+1}$), here is our algorithm.

Let $X$  be yet another variable.

{\bf Input}: An arbitrary Count Your Chickens! game, $G$, with $N+1$ locations, $K$ animals, and a given set
of blue squares, and positive integer $M$, and two variables $t$ and $X$.

{\bf Output}: A polynomial $P(X,t)$ of degree $M$ in $X$ and degree $N$ in $t$, such that
the coefficient of $X^i t^j$ is the probability of ending the game after exactly $i$ rounds with a capital of $j$ chicks.
It also outputs the probability of the game lasting longer than $M$ rounds.

We first initialize
$$
Q_0(X,t) := s^1=s \quad, \quad  R(X,t):=0  \quad ;
$$
and then for  $i=1 \dots, M$, we define, iteratively,
$$
Q'_i(X,t) := T( Q_{i-1}(X,t) )  \quad ;
$$
$$
Q_i(X,t) := Q'_i(X,t)- (Coefficient \quad of \quad s^{N+1} \quad in \quad  Q'_i(X,t))s^{N+1} \quad ;
$$
$$
R(X,t) := R(X,t)+(Coefficient \quad of \quad s^{N+1} \quad in \quad  Q'_i(X,t))X^i \quad .
$$

The output is $R(X,t)$, that tells us all the statistical information
for finishing in $\leq M$ rounds, and $\epsilon:=Q'_{M+1}(1,1)$ indicating the probability of {\bf not} terminating in $\leq M$ moves. 
You choose $M$ large enough so that $\epsilon$ is negligible.

Note that $R(X,t)$ is a polynomial of the two variables $X$ and $t$ of degrees $M$ and $N$ respectively.
This contains much more then {\it just} the probability of winning and the expected number of chicks.
If $\epsilon$ is tiny we can approximate the real thing by $R(X,t)$ and then $R(X,t)$ is the
bi-variate probability generating function of $(NumberOfRounds, NumberOfChicks)$ at the end.

(More precisely $R(X,t)/(1-\epsilon)$ is the {\it conditional}  probability generating function conditioned on
terminating in $\leq M$ rounds. From now on let $R(X,t):=R(X,t)/(1-\epsilon)$ ).

The probability of winning is the coefficient of $t^N$ in $R(1,t)$. The expected number of chicks is
$\frac{d}{dt}R(1,t)$. The expected number of rounds is $\frac{d}{dX}R(X,1)$. Similarly, we can
find the variances, the correlation, and any desired higher moments.

This is implemented in procedure {\tt ChSer(CB,t,X,M)} in the Maple package  {\tt CountChickens.txt}.
Typing

{\tt ChSer(CCb1(),t,X,60); }

gives the rather long $R(X,t)$ that can be seen in the output file

{\tt http://sites.math.rutgers.edu/\~{}zeilberg/tokhniot/oCountChickens1.txt} .

The more succinct command {\tt Info(CB,M)} uses $R(X,t)$ to extract the desired statistical information.

In particular the probability of winning turns out to be
$$
0.6410373996231\dots \quad,
$$
and the expected number of chicks, at the end of the game is
$$
39.32230439142343\dots  .
$$
If you have any doubts, we also have a {\it simulation} program that plays the game many times, and takes
the empirical averages. See

{\tt http://sites.math.rutgers.edu/~zeilberg/tokhniot/oCountChickens4.txt} \quad,

where the game is played one million times, and the empirical averages are very close to the above theoretical
values, confirming that the value of $39.22$ chicks in [MM] was a misprint.

We also found that the  variance of the number of chicks is
$1.2907513179745\dots$ that is rather small (explaining why the simulation values were so good),
the skewness is $-2.05489022\dots$ and the kurtosis is $7.8590953\dots$.

The average number of rounds happens to be $11.44706710\dots$ and its  variance is $6.28030112\dots$.
The correlation between the number of chicks and the number of rounds is                             
$-0.527785421907\dots$.

{\bf The more general Maple package UmbralMarkov.txt}

What we have here is what we call a {\it Weighted Markov Chain} with {\bf one} absorbing state.
A general weighted Markov chain with $n$ non-absorbing states and $s$ absorbing states is a directed graph
on $n+s$ vertices where the out-degree of each of the $s$ absorbing states is $0$, for each non-absorbing
state there is a probability distribution among its outgoing neighbors, and in addition each
edge carries a {\it weight}. You can think of the weight of a directed edge as the price that you have to pay every time you use it.
As you travel along this directed graph, according to the transition matrix, sooner or later you will wind
up in an absorbing state, and then you have a probability distribution regarding the total price of the travel,
for each of these absorbing states. You can also impose a ``minimum" total  price and a ``maximum" one,
like in the `Count Your Chickens!' game.

This more general scenario is implemented in the Maple package {\tt UmbralMarkov.txt} available from
the url  mentioned above, where there is also some sample output.

{\bf References}

[B] Peggy Brown, {\it ``Count Your Chickens!''}, a board game for ages $3+$, manufactured by Peaceable Kingdom.

[MM] David Mccune and Lori Mccune, {\it Counting your chickens with Markov chains}, Mathematics Magazine {\bf 92} (2019), 162-172.

\vfill\eject

\bigskip
\hrule
\bigskip
Shalosh B. Ekhad, c/o D. Zeilberger, Department of Mathematics, Rutgers University (New Brunswick), Hill Center-Busch Campus, 110 Frelinghuysen
Rd., Piscataway, NJ 08854-8019, USA. \hfill\break
Email: {\tt ShaloshBEkhad at gmail dot com}   \quad .
\bigskip
Doron Zeilberger, Department of Mathematics, Rutgers University (New Brunswick), Hill Center-Busch Campus, 110 Frelinghuysen
Rd., Piscataway, NJ 08854-8019, USA. \hfill\break
Email: {\tt DoronZeil at gmail  dot com}   \quad .
\bigskip
\hrule
\bigskip
{\bf Exclusively published in the Personal Journal of Shalosh B.  Ekhad and Doron Zeilberger and arxiv.org \quad .}
\bigskip
\hrule
\bigskip
Written: July 18, 2019.
\end